\PassOptionsToPackage{table}{xcolor}
\documentclass[AMA,twocolumn]{USG}
\usepackage{graphicx} 
\graphicspath{{./images/}}
\usepackage{amsmath}
\usepackage{amssymb}
\usepackage{mathtools}
\usepackage{algorithm}
\usepackage{algpseudocode}
\usepackage{multirow}
\usepackage[short]{optidef}
\usepackage{lettersp}
\usepackage{optidef}
\usepackage{nicematrix}
\usepackage{bm}
\usepackage{comment}
\usepackage{tikz}
\usetikzlibrary{arrows.meta, positioning, shapes.geometric, calc}

\articletype{Short Communication}%

\received{Date Month Year}
\revised{Date Month Year}
\accepted{Date Month Year}
\journal{Journal}
\volume{00}
\copyyear{2023}
\startpage{1}
\articledoi{10.1002/asna.20230101} 

\begin{document}

\title{Exploiting Micro-Structure in the Generalized Riccati Recursion for Optimal Control}
\author[1,2]{Louis Callens}
\author[1,2]{Wilm Decré}
\address[1]{\orgdiv{MECO Research team, Department of Mechanical Engineering}, \orgname{KU Leuven}, \orgaddress{\state{Leuven}, \country{Belgium}}}
\address[2]{\orgdiv{Flanders Make @ KU Leuven}, \orgaddress{\state{Leuven}, \country{Belgium}}}

\corres{Corresponding author Louis Callens, \email{louis.callens@kuleuven.be}}
\presentaddress{This is sample for present address text this is sample for present address text.}

\fundingInfo{This work has been carried out within the framework of the Flanders Make SBO project LearnOptra. Flanders Make is the Flemish strategic research centre for the manufacturing industry.}

\abstract[Abstract]{Optimal Control Problems (OCPs) have been applied to a variety of applications. Such problems are often transcribed into a Nonlinear Program (NLP) for which solution methods typically solve a Linear Quadratic (LQ) subproblem. To efficiently solve this LQ subproblem, the Riccati recursion exploits the block-diagonal structure that arises from discretizing an OCP using a Multiple Shooting method. However, reformulating general dynamics constraints to fit the required problem structure to use the Riccati recursion introduces additional micro-structure within these block matrices that remains unexploited. This work demonstrates how structure-exploiting solvers such as Fatrop can benefit from exploiting not only the block-diagonal structure but also the micro-structure within those blocks. We propose a structure-exploiting LU decomposition algorithm and integrate it within the Riccati recursion. Numerical results on randomized linear systems and two OCP examples demonstrate a reduction in computation time of the search direction of up to 20\%, with the largest speedup observed for relative zero block sizes around 0.25 and large stage-wise equality constraint jacobians.}



\maketitle

\section{Introduction}
Optimal Control has been widely used in various applications such as robotics. In these applications, robots solve the Optimal Control Problem (OCP) to find control inputs that lead to a feasible and optimal trajectory. The computation time required to solve the OCP is of interest for a number of reasons. Firstly, to enable the robot to react to changes in the environment or to deal with plant-model mismatch, the OCP is repeatedly solved in a receding horizon fashion, which leads to (Nonlinear) Model Predictive Control ((N)MPC) and imposes real-time constraints on the computation time since the OCP must be solved at a fixed rate. Secondly, OCPs can be used offline to generate training data\cite{campd}. Even though there is no real-time constraint on the computation time for offline data generation, computational efficiency is still of interest due to the large amount of data required. Thirdly, high fidelity simulations to gain insight in specific applications or dynamics can be computationally expensive.

Direct methods discretize continuous-time OCPs to obtain a Nonlinear Program (NLP). 

Two popular methods to solve NLPs with both equality and inequality constraints are interior-point methods and sequential quadratic programming (SQP) methods. The former, such as IPOPT\cite{ipopt}, solve a related NLP, called the barrier problem and drive the barrier parameter gradually to zero to obtain a solution to the original NLP. The system describing the Karash-Kuhn-Tucker (KKT) conditions is solved to find the optimal solution. The latter, such as qpOASES\cite{qpoases}, solve a series of Quadratic Programs (QP) that each approximate the original problem around the current iterate. The QP itself is solved either using an active-set method, which attempts to find the constraints that are active and considers those as equality constraints or using an interior-point method.

Either way, a Linear Quadratic (LQ) problem is solved as a subproblem which is often the most time-consuming step of the solver. For OCPs, the KKT-matrix of the LQ problem has a block-diagonal structure. This structure arises from a multiple shooting transcription of the OCP, where the states are kept as optimization variables and dynamics constraints are added which typically converges faster for nonlinear problems than eliminating the states as in a single shooting approach. Additionally, stage-wise objectives and constraints are typically considered, which leads to the block-diagonal structure of the KKT matrix that can be exploited using the Riccati recursion \cite{ricatti_rao, ricatti-frison}. HPIPM \cite{hpipm} implements the Riccati recursion to solve QPs efficiently where BLASFEO \cite{blasfeo} is used to perform linear algebra operations efficiently. However, HPIPM cannot directly deal with stage-wise equality constraints. Therefore, a generalization of the Riccati recursion which can deal with such constraints directly has been proposed \cite{fatrop-ricatti}. Fatrop implements this generalization of the recursion efficiently using BLASFEO and integrates it within an interior-point solver\cite{fatrop}.

Fatrop\cite{fatrop} solves the Constrainted Optimal Control Problem (COCP) given by
\begin{mini!}{x_k, u_k, x_K}{l_K(x_K) + \sum_{k=0}^{K-1}{l_k(u_k, x_k)}}{\label{eq:fatrop-ocp}}{}
    \addConstraint{}{x_{k+1} = f_k(u_k, x_k)} \label{eq:fatrop-ocp-dynamics}
    \addConstraint{}{L_k \leq g_k(u_k, x_k) \leq U_k}
    \addConstraint{}{L_K \leq g_K(x_K) \leq U_K}
    \addConstraint{}{h_k(u_k, x_k) = 0}
    \addConstraint{}{h_K(x_K) = 0}
\end{mini!}
for $k = 0, 1, \ldots K-1$ with $K$ the horizon length, $x_k$ the state variables, $u_k$ the control variables. The objective is stage-wise, as well as the equality and inequality constraints. The dynamics constraint \eqref{eq:fatrop-ocp-dynamics} couples stage $k$ with stage $k+1$.

Note that more general dynamics constraints, where $x_{k+1}$ appears implicitly, do not fit the required problem structure. However, such dynamics constraints arise when using collocation methods for transcription or when implicit integrators, known for higher numerical stability, are used. 


\subsection{Contributions}
In this work, we show how reformulating general dynamics constraints to fit the required problem structure introduces additional micro-structure within these block matrices that remains unexploited. We propose a structure-exploiting LU decomposition algorithm that is integrated within the Riccati recursion and validate the effectiveness of this approach on randomized KKT-systems and two OCP examples.

\subsection{Structure of the paper}
The remainder of this paper is structured as follows. Section \ref{sec:recursion} presents the Riccati recursion and specifically considers a matrix decomposition as a key step in the algorithm. Section \ref{sec:reformulation} shows how general dynamics can be  can be reformulation to fit the required structure. Section \ref{sec:structure-Gu} analyzes the structure in the equality constraint jacobian introduced by this reformulation. Section \ref{sec:structure-exploiting-lu} presents the structure-exploiting LU-decomposition. Numerical results are demonstrated in Section \ref{sec:results} and conclusions are drawn in Section \ref{sec:conclusion}.
\section{LU decomposition in the Riccati recursion} \label{sec:ricatti}
The generalized Riccati recursion solves the KKT-system of the problem defined in \eqref{eq:fatrop-ocp} while exploiting the structure of that KKT-system. The structure for horizon lenght $K = 2$ is given by\cite{fatrop}
\begin{equation}
    \begin{bNiceArray}[first-row]{cccccccccc:c}
        x_2     & \nu_2 & \pi_2 & u_1 & x_1 & \lambda_1 & \pi_1 & u_0 & x_0 & \lambda_0 & \\
        Q_2     & G_2  & -\bm{I} &          &          &          &         &          &          &          & q_2 \\
        G_2     &      &         &          &          &          &         &          &          &          & g_2 \\
        -\bm{I} &      &         & B_1      & A_1      &          &         &          &          &          & b_1 \\
                &      & B_1'    & R_1      & S_1'     & G_{1,u}' &         &          &          &          & r_1 \\
                &      & A_1'    & S_1      & Q_1      & G_{1,x}' & -\bm{I} &          &          &          & q_1 \\
                &      &         & G_{1,u}  & G_{1,x}  &          &         &          &          &          & g_1 \\
                &      &         &          & -\bm{I}  &          &         & B_0      & A_0      &          & b_0 \\
                &      &         &          &          &          & B_0'    & R_0      & S_0'     & G_{0,u}' & r_0 \\
                &      &         &          &          &          & A_0'    & S_0      & Q_0      & G_{0,x}' & q_0 \\
                &      &         &          &          &          &         & G_{0,u}  & G_{0,x}  &          & g_0 \\
    \end{bNiceArray}.
    \label{eq:KKT}
\end{equation}
The Riccati recursion considers each block to be dense and solves the system using a backward recursion algorithm, followed by a forward substution. The backward recursion iteratively eliminates $x_m$, $\nu_m$, $\pi_m$ and $u_{m-1}$ where $m$ is the largest remaining stage-index. Using the first and third block row of the remaining KKT-system, $x_m$ and $\pi_m$ are eliminated. Constraints on $x_m$ are transferred to stage $m-1$ using the dynamics equation which leads to new stagewise constraint jacobians $\overline{G}_{u,m-1}$ and $\overline{G}_{x,m-1}$. The former is then decomposed as
\begin{equation}
    \overline{G}_{u,m-1} = T_{m-1,L} \begin{bmatrix} -\bm{I}_\rho & \\ & \bm{0} \end{bmatrix} T_{m-1,R} \label{eq:decomposition}
\end{equation}
with $T_{m-1,L}$ and $T_{m-1,R}$ invertible and $\rho$ the rank of $\overline{G}_{u,m-1}$. This decomposition is used to eliminate $\rho$ controls (the remaining controls are eliminated using a Shur complement step) and can be obtained after performing an LU-decomposition of $\overline{G}_{u,m-1}$. 

\begin{figure}
    \centering
    \includegraphics[width=\linewidth]{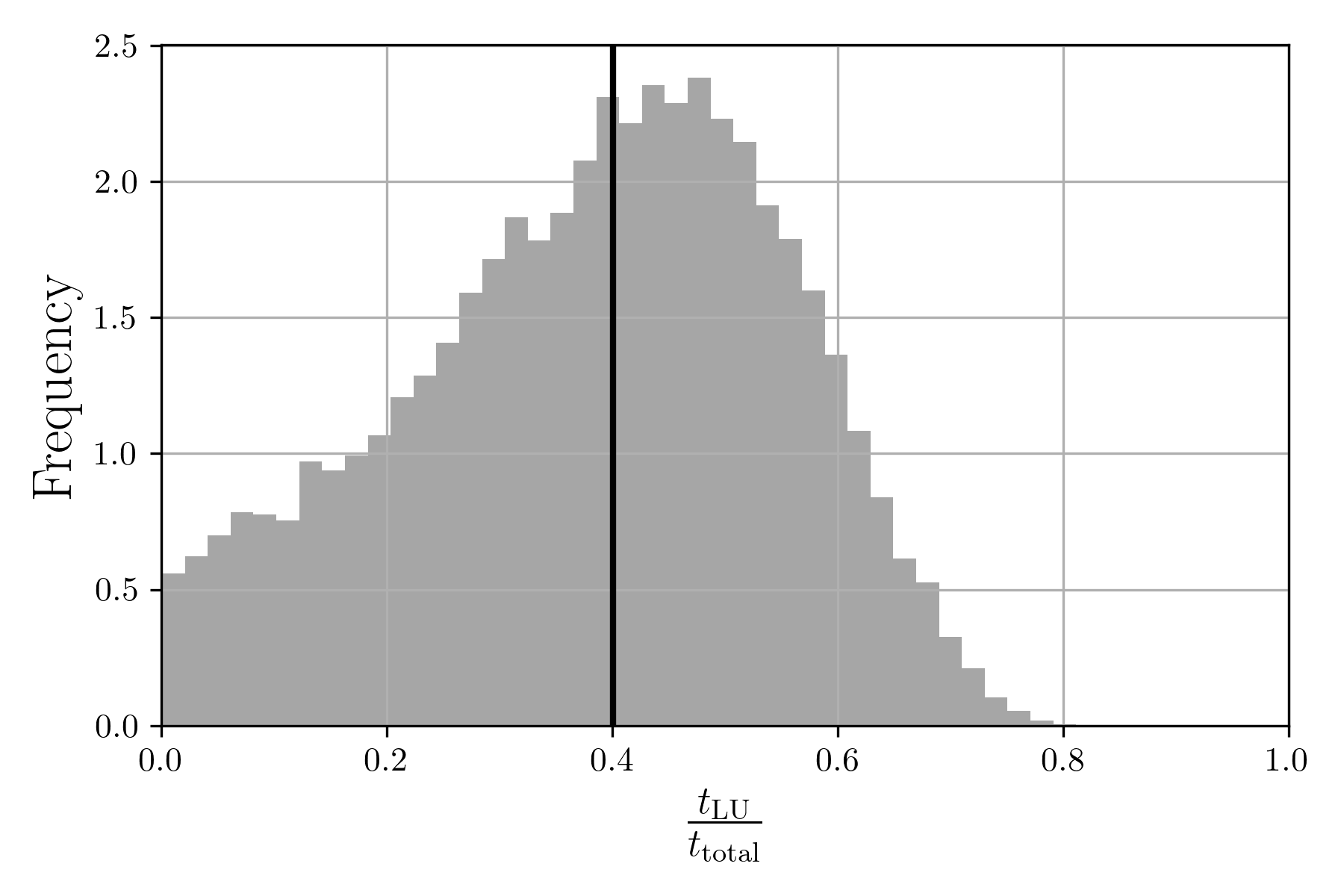}
    \caption{Distribution of the fraction of the computation time of the Riccati recursion algorithm spent performing LU decompositions for the linear systems from the benchmark systems.}
    \label{fig:lu_relevance_distribution}
\end{figure}

Figure \ref{fig:lu_relevance_distribution} shows a distribution of the fraction of the computation time of the Riccati recursion algorithm spent performing LU decompositions for the linear systems from the benchmark systems described in Section \ref{sec:recursion}. The median fraction is around 40\%, which motivates the need for a more efficient decomposition algorithm. 

\section{Methodology} \label{sec:reformulation}
A generalization of the dynamics constraint \eqref{eq:fatrop-ocp-dynamics} is given by
\begin{equation}
    f_k(u_k, x_k, x_{k+1}) = 0 \label{eq:general-dynamics}
\end{equation}
where $f_k : \mathbb{R}^{n_u^k} \times \mathbb{R}^{n_x^k} \times \mathbb{R}^{n_x^{k+1}} \to \mathbb{R}^{n_f^k}$. After introducing auxiliary controls $z_k$, these dynamics can be rewritten as
\begin{equation}
    x_{k+1} = z_k
    \quad \text{ s.t. } \quad
    f_k(u_k, x_k, z_{k}) = 0 \label{eq:additional-constraints}.
\end{equation}
This reformulated problem fits the required problem formulation at the cost of introducing additional variables and stage-wise equality constraints. 

More generally, suppose the states can be written as 
\begin{equation}
    x_{k+1} = \begin{bmatrix} z_k \\ y_k \end{bmatrix}
\end{equation} with $z_k \in \mathbb{R}^{n_z^k}$ representing states that do not fit the standard form in \eqref{eq:fatrop-ocp-dynamics} and $y_k \in \mathbb{R}^{n_x^{k+1}-n_z^k}$ representing the states that can be eliminated explicitly. The dynamics can then be written as
\begin{equation}
    f_k(u_k, x_k, x_{k+1}) = 
    \begin{bmatrix}
        f_{k,1}\left(u_k, x_k\right) \\
        f_{k,2}\left(u_k, x_k, z_k\right) \\
        y_{k} - f_{k,3}\left(u_k, x_k, z_k\right)
    \end{bmatrix}
    = 0 \label{eq:dynamics-decomposed}
\end{equation}
where $f_{k,1} : \mathbb{R}^{n_u^k} \times \mathbb{R}^{n_x^k} \to \mathbb{R}^{n_{f_1}^k}$, $f_{k,2} : \mathbb{R}^{n_u^k} \times \mathbb{R}^{n_x^k} \times \mathbb{R}^{n_z^k} \to \mathbb{R}^{n_{f_2}^k}$ and $f_{k,3} : \mathbb{R}^{n_u^k} \times \mathbb{R}^{n_x^k} \times \mathbb{R}^{n_z^k} \to \mathbb{R}^{n_{f_3}^k}$ where $n_{f_1}^k + n_{f_2}^k + n_{f_3}^k = n_f^k$ and $n_{f_3}^k = n_x^{k+1} - n_z^k$. The dynamics equation \eqref{eq:dynamics-decomposed} is equivalent to
\begin{subequations}
    \label{eq:reformulated-dynamics}
\begin{align}
    x_{k+1} &= 
    \begin{bmatrix}
        z_k \\
        f_{k,3}\left(u_k, x_k, z_k\right)
    \end{bmatrix}
    \coloneqq \tilde{f}_k(u_k, x_k, z_k)\\
    &\text{ s.t. } \quad  \begin{bmatrix} f_{k,1}\left(u_k, x_k\right)\\ f_{k,2}\left(u_k, x_k, z_k\right) \end{bmatrix} = 0
\end{align}
\end{subequations}
such that the problem structure fits that of \eqref{eq:fatrop-ocp}. Note that variables $y_k$ are not actually introduced as optimization variables but variables $z_k$ are introduced as controls.

\subsection{Structure of $\overline{G}_u$} \label{sec:structure-Gu}
If the reformulation shown in \eqref{eq:reformulated-dynamics} is used, the matrix $\overline{G}_u$ is given by
\begin{align}
    \overline{G}_u = 
    \begin{bmatrix}
        \nabla_{u_k} g & \nabla_{z_k} g \\
    \end{bmatrix}
    \coloneqq
    \begin{bmatrix}
        M_1 &  \\
        M_2 & M_3 \\
        M_4 & M_5
    \end{bmatrix}
    \label{eq:Gu-structure}\\
    \intertext{where}
    g = \begin{bmatrix} h_k(u_k, x_k) \\ f_{k,1}(u_k, x_k) \\ f_{k,2}(u_k, x_k, z_k)\\ H_{k+1} \tilde{f}_k(u_k, x_k, z_k) \end{bmatrix}
\end{align}
where $M_1 \in \mathbb{R}^{n_h^k + n_{f_1}^k \times n_u^k}$, $M_2 \in \mathbb{R}^{n_{f_2}^k \times n_u^k}$, $M_3 \in \mathbb{R}^{n_{f_2}^k \times n_z^k}$, $M_4 \in \mathbb{R}^{n_c^k \times n_u^k}$ and $M_5 \in \mathbb{R}^{n_c^k \times n_z^k}$ where $n_c^k$ represents the number of constraints carried over the next stage and is equal to the number of rows in $H_{k+1}$. Note that no assumptions are made about invertability or relative dimensions of any of the matrices $M_i$.

The size of the top-right block is given by $(n_h^k + n_{f_1}^k) n_z^k$, while the size of the total matrix is given by $S_\mathrm{total} \coloneqq (n_u^k + n_z^k) (n_h^k + n_{f_1}^k + n_{f_2}^k + n_c^k)$, meaning the relative size of the zero block compared to the full matrix is given by 
\begin{equation}
S_\mathrm{rel} \coloneqq ((n_h^k + n_{f_1}^k) n_z^k)/S_\mathrm{total}. \label{eq:Srel}
\end{equation}

If $n_z^{k} = n_x^{k+1}$, meaning no elements from $x_{k+1}$ enter $f_k$ explicitly, then $\tilde{f}_k$ as defined in \eqref{eq:reformulated-dynamics} does not depent on $u_k$ and hence $M_4 = 0$. 

The decomposition shown in \eqref{eq:decomposition} can be obtained by performing an LU-decomposition first\cite{fatrop-ricatti}. The structure in $\overline{G}_u$ can be exploited by performing a structure-exploiting LU-decomposition.

\subsection{Structure-Exploiting LU-decomposition} \label{sec:structure-exploiting-lu}

Instead of decomposing the matrix in \eqref{eq:Gu-structure} directly, an LU-decomposition with full pivoting of $M_1$ only can be performed first, leading to $P_1 M_1 Q_1 = L_1 U_1$ where
\begin{align}
    L_1 &= \begin{bmatrix} K_1 & \\ K_2 & \bm{I}_{n_g-\rho_1} \end{bmatrix} & 
    &\text{and} &
    U_1 &= \begin{bmatrix} V_1 & V_2 \\ & \bm{0}_{n_g-\rho_1 \times n_u-\rho_1} \end{bmatrix}
\end{align}
where $\rho_1$ is the rank of $M_1$. Because only pivots from within $M_1$ are selected, the structure is preserved and we obtain
\begin{equation}
    \begin{bmatrix} P_1 & \\ & \bm{I}_{n_f+n_c} \end{bmatrix} 
    \overline{G}_u
    \begin{bmatrix} Q_1 & \\ & \bm{I}_{n_z} \end{bmatrix} 
    = \begin{bmatrix}
        L_1 & \\
         & \bm{I}_{n_f+n_c}
    \end{bmatrix}
    \begin{bmatrix}
        U_1 & \\
        M_2 Q_1 & M_3\\ 
        M_4 Q_1 & M_5
    \end{bmatrix}.
\end{equation}
Define
\begin{align}
    \begin{bmatrix} K_3 & K_5\\ K_4 & K_6 \end{bmatrix} \coloneqq \begin{bmatrix} M_2 \\ M_4 \end{bmatrix} Q_1 \begin{bmatrix} V_1^{-1} & -V_1^{-1} V_2 \\ & \bm{I}_{n_u - \rho_1} \end{bmatrix},
\end{align}
such that
\begin{align}
    \begin{bmatrix} K_3 & \bm{0}\\ K_4 & \bm{0} \end{bmatrix} U_1 + \begin{bmatrix} \bm{0} & K_5 \\ \bm{0} & K_6 \end{bmatrix} = \begin{bmatrix} M_2 \\ M_4 \end{bmatrix} Q_1
\end{align}
holds. Hence, we can see that
\begin{equation}
    \begin{bmatrix} P_1 & \\ & \bm{I}_{n_{f}+n_c} \end{bmatrix} 
    \overline{G}_u
    \begin{bmatrix} Q_1 & \\ & \bm{I}_{n_z} \end{bmatrix} 
    = \begin{bmatrix}
        L_1 & \\
        \begin{bmatrix} K_3 & \bm{0} \\ K_5 & \bm{0} \end{bmatrix} & \bm{I}_{n_{f}+n_c}
    \end{bmatrix}
    \begin{bmatrix}
        U_1 &\\
        \begin{bmatrix} \bm{0} & K_5 \end{bmatrix} & M_3\\
        \begin{bmatrix} \bm{0} & K_6 \end{bmatrix} & M_5
    \end{bmatrix}.
\end{equation}
After a row permutation, only the decomposition
\begin{equation}
    P_2\begin{bmatrix} K_5 & M_3\\ K_6 & M_5\end{bmatrix}Q_2 = L_2 U_2
\end{equation}
is required to obtain the full decomposition of $\overline{G}_u$. Algorithm \ref{alg:structure-exploiting-lu-carry} summarizes the procedure to obtain the decomposition. Matrix dimensions are shown in Table \ref{tab:matrix-dimensions}. Appendix \ref{app:structure-exploiting-lu} contains a validation of the algorithm.

\begin{algorithm}
    \caption{Structure-Exploiting LU decomposition of matrix $\overline{G}_u$}
    \label{alg:structure-exploiting-lu-carry}
    \begin{algorithmic}
    \Require Matrix $\overline{G}_u$ given by \eqref{eq:Gu-structure}
    \Ensure LU with complete pivoting: $P \overline{G}_u Q = L U$
    \State Compute LU decomposition $P_1 M_1 Q_1 = \begin{bmatrix}K_1 & \\ K_2 & \bm{I}\end{bmatrix} \begin{bmatrix} V_1 & V_2\\ & \bm{0}\end{bmatrix}$
    \State $\begin{bmatrix} K_3 & K_5\\ K_4 & K_6 \end{bmatrix} \gets \begin{bmatrix} M_2 \\ M_4 \end{bmatrix} Q_1 \begin{bmatrix} V_1^{-1} & -V_1^{-1} V_2 \\ & \bm{I}_{n_u - \rho_1} \end{bmatrix}$
    \State Compute LU decomposition $P_2 \begin{bmatrix} K_5 & M_3\\ K_6 & M_5 \end{bmatrix} Q_2 = L_2 U_2$
    \State $K_7 \gets P_2 \begin{bmatrix} K_3\\ K_4 \end{bmatrix}$
    \State $\tilde{V}_2 \gets \begin{bmatrix} V_2 & 0 \end{bmatrix} Q_2$
    \State $L \gets \begin{bmatrix} K_1 & & \\ K_7 & L_2 & \\ K_2 &  & \bm{I}_{n_g-\rho_1} \end{bmatrix}$
    \State $U \gets \begin{bmatrix} V_1 & \tilde{V}_2\\ & U_2\\ & \end{bmatrix}$
    \State $P \gets \begin{bmatrix} \bm{I}_{\rho_1} & & \\ & P_2 & \\ & & \bm{I}_{n_g-\rho_1} \end{bmatrix} \begin{bmatrix} \bm{I}_{\rho_1} & & \\ & & \bm{I}_{n_f+n_c} \\ & \bm{I}_{n_g-\rho_1} & \end{bmatrix} \begin{bmatrix} P_1 & \\ & \bm{I}_{n_f+n_c} \end{bmatrix}$
    \State $Q \gets \begin{bmatrix} Q_1 & \\ & \bm{I}_{n_z}\end{bmatrix} \begin{bmatrix}\bm{I}_{\rho_1} &\\ & Q_2 \end{bmatrix}$

    \end{algorithmic}
\end{algorithm}

\begin{table}
    \centering
    \rowcolors{2}{gray!30}{white}
    \caption{Dimensions of matrices appearing in algorithm \ref{alg:structure-exploiting-lu-carry}}
    \label{tab:matrix-dimensions}
    \begin{tabular}[width=\linewidth]{ccc}
        \rowcolor{gray!70}
         & \textbf{Rows} & \textbf{Columns} \\
        $M_1$ & $n_h + n_{f_1}$ & $n_u$ \\
        $M_2$ & $n_{f_2}$ & $n_u$ \\
        $M_3$ & $n_{f_2}$ & $n_z$ \\
        $M_4$ & $n_c$ & $n_u$ \\
        $M_5$ & $n_c$ & $n_z$ \\
        $\overline{G}_u$ & $n_h + n_{f_1} + n_{f_2} + n_c$ & $n_u + n_z$ \\
        $P_1$ & $n_h + n_{f_1}$ & $n_h + n_{f_1}$ \\
        $Q_1$ & $n_u$ & $n_u$ \\
        $K_1$ & $\rho_1$ & $\rho_1$ \\
        $K_2$ & $n_g - \rho_1$ & $\rho_1$ \\
        $V_1$ & $\rho_1$ & $\rho_1$ \\
        $V_2$ & $\rho_1$ & $n_u - \rho_1$ \\
        $L_1$ & $n_h + n_{f_1}$ & $n_h + n_{f_1}$ \\
        $U_1$ & $n_h + n_{f_1}$ & $n_u$ \\
        $K_3$ & $n_{f_2}$ & $\rho_1$ \\
        $K_4$ & $n_c$ & $\rho_1$ \\
        $K_5$ & $n_{f_2}$ & $n_u - \rho_1$ \\
        $K_6$ & $n_c$ & $n_u - \rho_1$ \\
        $P_2$ & $n_{f_2} + n_c$ & $n_{f_2} + n_c$ \\
        $Q_2$ & $n_u - \rho_1 + n_z$ & $n_u - \rho_1 + n_z$ \\
        $L_2$ & $n_{f_2} + n_c$ & $n_{f_2} + n_c$ \\
        $U_2$ & $n_{f_2} + n_c$ & $n_u - \rho_1 + n_z$ \\
        $K_7$ & $n_{f_2} + n_c$ & $\rho_1$ \\
        $\tilde{V}_2$ & $\rho_1$ & $n_u - \rho_1 + n_z$
    \end{tabular}
\end{table}
    
\section{Results and discussion} \label{sec:results}
The proposed method is evaluated extensively on randomized linear systems, as explained in section \ref{sec:recursion}. Additionally, the method is applied to OCP examples in section \ref{sec:ocp_example}.

\subsection{Effect on Riccati Recursion} \label{sec:recursion}
To evaluate the effect of the structure-exploiting LU decomposition on the Riccati recursion randomized linear systems of the form shown in \eqref{eq:KKT} are constructed. 
We uniformly sample $K$ from $U(2, 12)$, $n_u$, $n_h$, $n_z$, $n_v$, $n_g$ from $U(0, 30)$, and $n_x$ from $U(n_z, n_z + 30)$. Block matrices are generated randomly as well. A total of 20000 linear systems are generated.

\begin{figure}
    \centering
    \includegraphics[width=\linewidth]{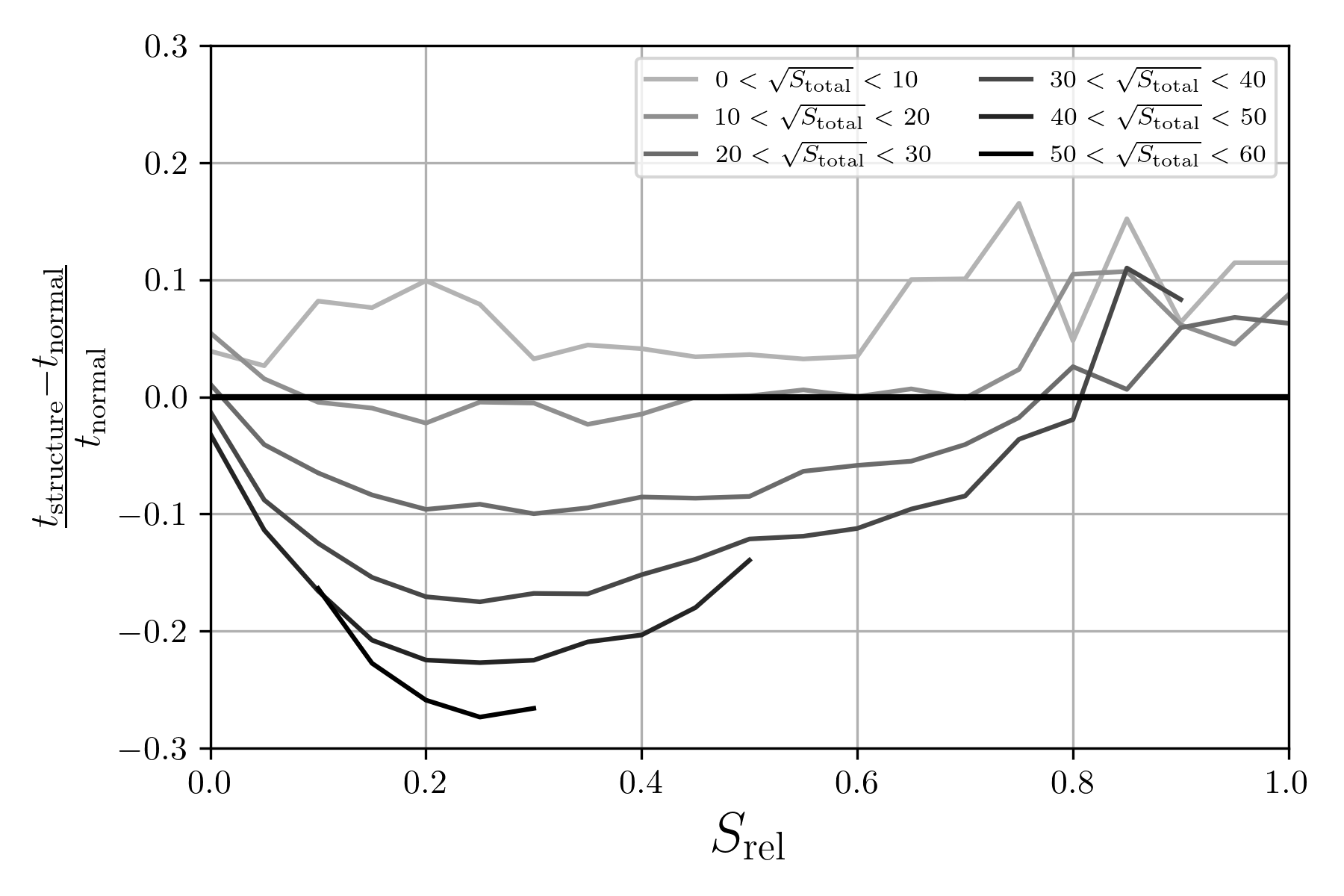}
    \caption{Relative improvement of the structure-exploiting LU decomposition compared to the normal LU decomposition as a function of $S_\mathrm{rel}$. The computation times $t_\mathrm{normal}$ and $t_\mathrm{structure}$ represent the total computation time of the Riccati recursion using the normal and structure-exploiting LU decomposition respectively. Mean values are shown for different sizes of the full matrix.}
    \label{fig:relative_improvement_relative_area}
\end{figure}

Figure \ref{fig:relative_improvement_relative_area} shows the relative difference in computation time of the Riccati recursion between the structure-exploiting LU decomposition and the standard LU decomposition. These results are shown as a function of $S_\mathrm{rel}$, the relative size of the top-right zero block defined in \eqref{eq:Srel}. Mean values are shown for different binned sizes of the full matrix. As the size of the full matrix increase, the speedup of the structure-exploiting LU decomposition increases as well. The speedup is largest for a relative size of around 0.25. For smaller relative sizes, there is less structure to be exploited. For higher relative sizes, the rank of $\overline{G}_{u,m-1}$ decreases and the LU decomposition can be terminated early which leads to a smaller fraction of the overal time spent in the LU decomposition.

\subsection{OCP Examples} \label{sec:ocp_example}
Two OCPs are considered to evaluate the effect of the structure-exploiting LU decomposition. These are detailed in Section \ref{sec:ocp_description}. Results are discussed in Section \ref{sec:ocp_results}.
\subsubsection{Description of the examples} \label{sec:ocp_description}
The first OCP is a chain of $n$ pendulums. The dynamics are modelled as a Differential Algebraic Equation (DAE)\cite{dea-multibody} given by
\begin{align}
    M(q) \ddot{q} &= f(q, F) - J_c(q)^\top w\\
    0 &= c(q) \label{eq:pendulum-constraint}
\end{align}
where $q \in \mathbb{R}^{3n}$ are the positions of the masses, $M(q) \coloneqq 1$kg$\bm{I}_{3n}$ is the mass matrix, $f(q, F)$ is the vector of forces acting on the masses (including gravity and external forces), $c(q)$ enforces the distance between two consequtive masses to be equal to $L = 0.5$m, $J_c(q)$ is the Jacobian of the constraints and $w \in \mathbb{R}^n$ is a vector of Lagrange multipliers. To obtain an index-1 DAE, the second derivative of \eqref{eq:pendulum-constraint} is enforced to be zero, which leads to the equations
\begin{align}
    M(q) \ddot{q} &= f(q, \dot{q}, F) - J_c(q)^\top w\\
    0 &= J_c(q) \ddot{q} +  \frac{\partial}{\partial q} \left(J_c(q)^\top \dot{q}\right)^\top\dot{q}
\end{align}
From the first equation, we can compute $\ddot{q}$ as a function of $q$, $\dot{q}$, $w$, and $F$. The second equation enforces the second derivative of the lengths of the pendulums (given by the distances between two consecutive masses) to be constant. To prevent constraint drift, Baumgartne stabilization \cite{baumgarte_stabilization_1972} is used which adds stiffness to the system. A symplectic integrator is used to discretize the dynamics leading to the formulation
\begin{align}
    \begin{bmatrix} \dot{q}_{k+1}\\ q_{k+1}\end{bmatrix} &= \begin{bmatrix} \dot{q}_k + \ddot{q}_{k}(q_k, \dot{q}_k, w_k, F) \Delta t\\ q_k + \dot{q}_{k+1} \Delta t \end{bmatrix}\\\Longleftrightarrow \begin{bmatrix} \dot{q}_{k+1}\\ q_{k+1}\end{bmatrix} &= \begin{bmatrix} z_k\\ q_k + z_k \Delta t \end{bmatrix} \text{ s.t. } z_k = \dot{q}_{k} + \ddot{q}_{k}(q_k, \dot{q}_k, w_k, F)
\end{align}
where both formulations are equivalent, but the second one fits the required structure by introducing variables $z_k \in \mathbb{R}^{3n}$ and adding equality constraints.

The objective of the OCP is to minimize terminal velocity and deviation from the configuration at rest at the final time, starting from a given initial configuration, while minimizing control effort throughout the trajectory. Only the third and fifth masses can be actuated by applying a force in 3D. This leads to $6 + n$ controls (from $F$ and $w$).

Table \ref{tab:block_dimensions} shows the dimensions of the problem as a function of the number of pendulums $n$. This means the relative size of the top-right zero block is given by $\frac{3n^2}{16n^2+24n}$, which for $n=8$ is equal to 0.158.

The second OCP example is a simulation of a human periodic walking gait where a muscluloskeletal model is used \cite{biomechanical} \footnote{The dynamics can be constructed using the open-source software \url{https://codeberg.org/Lars-DHondt/SOC_walking_DHondt2026} }. A simplified model of only the lower legs is used. The joint configuration at time step $k$ is expressed by $q_k \in \mathbb{R}^9$. To enforce periodicity using only stage-wise equality constraints, the initial joint configuration and joint velocities have to be propagated through the horizon. This leads to a state vector given by $x_k = \begin{bmatrix} q_k^\top & \dot{q}_k^\top & \psi_k^\top \end{bmatrix}^\top$ where $\psi_k$ are auxiliary variables to propagate the initial state. The controls consist of muscle activation $a_k \in \mathbb{R}^{18}$. A third-order collocation scheme is used to integrate the dynamics leading to additional controls $\bm{c}_{k} \coloneqq \begin{bmatrix} c_{k,1}^\top & c^\top_{k,2} & c^\top_{k,3} \end{bmatrix}^\top \in \mathbb{R}^{27}$ representing the joint configuration at the collocation points and the derivatives $\bm{\dot{c}}_k$ and $\bm{\ddot{c}}_k$.
Specifically, the dynamics and corresponding equations are enforced using
\begin{align}
    \begin{bmatrix} 
        q_{k+1}\\ \dot{q}_{k+1}\\ \psi_{k+1}
    \end{bmatrix}
    &= 
    \begin{bmatrix} 
        c_{k,3}\\ \dot{c}_{k,3}\\ \psi_k
    \end{bmatrix}
    \text{ s.t. }
    \begin{bmatrix} 
        h_\mathrm{coll}(q_k, \bm{c}_k) - \bm{\dot{c}}_k\\
        h_\mathrm{coll}(\dot{q}_k, \bm{\dot{c}}_k) - \bm{\ddot{c}}_k\\
        f_\mathrm{dyn}(q_k, \dot{q}_k, \bm{c}_k, \bm{\dot{c}}_k, \bm{\ddot{c}}_k, a_k)
    \end{bmatrix}
    = 0. \label{eq:walking_gait_dynamics}
\end{align}

Comparing \eqref{eq:walking_gait_dynamics} with \eqref{eq:reformulated-dynamics} after arranging some of the collocation constraints, we can see that $z_k = \begin{bmatrix} c_{k,3}^\top & \dot{c}_{k,3}^\top \end{bmatrix}^\top$ and $u_k = \begin{bmatrix} a_k^\top & c_{k,1}^\top & c_{k,2}^\top & \dot{c}_{k,1}^\top & \dot{c}_{k,2}^\top & \bm{\ddot{c}}_{k}^\top \end{bmatrix}^\top$. The relative size of the top-right zero block is 0.136. Table \ref{tab:block_dimensions} shows relevant dimensions of this problem.

\begin{table}
    \centering
    \rowcolors{2}{gray!30}{white}
    \caption{Dimensions of the blocks in $\bar{G}_u$ for the OCP example.}
    \label{tab:block_dimensions}
    \begin{tabular}{ccc}
        \rowcolor{gray!70}
        \textbf{Dimension} & \textbf{$n$-Pendulum} & \textbf{Walking gait} \\
        $n_h$ & $n$ & 0 \\
        $n_{f_1}$ & 0 & 4 $n_q$ \\
        $n_{f_2}$ & $3n$ & $36$ \\
        $n_{f_3}$ & $3n$ & $2n_q-1$ \\
        $n_u$ & $n + 6$ & $2n_q + 7n_q$ \\
        $n_z$ & $3n$ & $2n_q$ \\
        $n_x$ & $6n$ & $2n_q + 2n_q - 1$\\
        $S_\mathrm{rel}$ & $0.158$ (for $n=8$) & $0.136$ \\
        $\sqrt{S_\mathrm{total}}$ & $34.87$ (for $n=8$) & 76.37
    \end{tabular}
\end{table}

\subsubsection{OCP Example Results} \label{sec:ocp_results}
\begin{figure}
    \centering
    \subfloat{
    \includegraphics[width=0.8\linewidth]{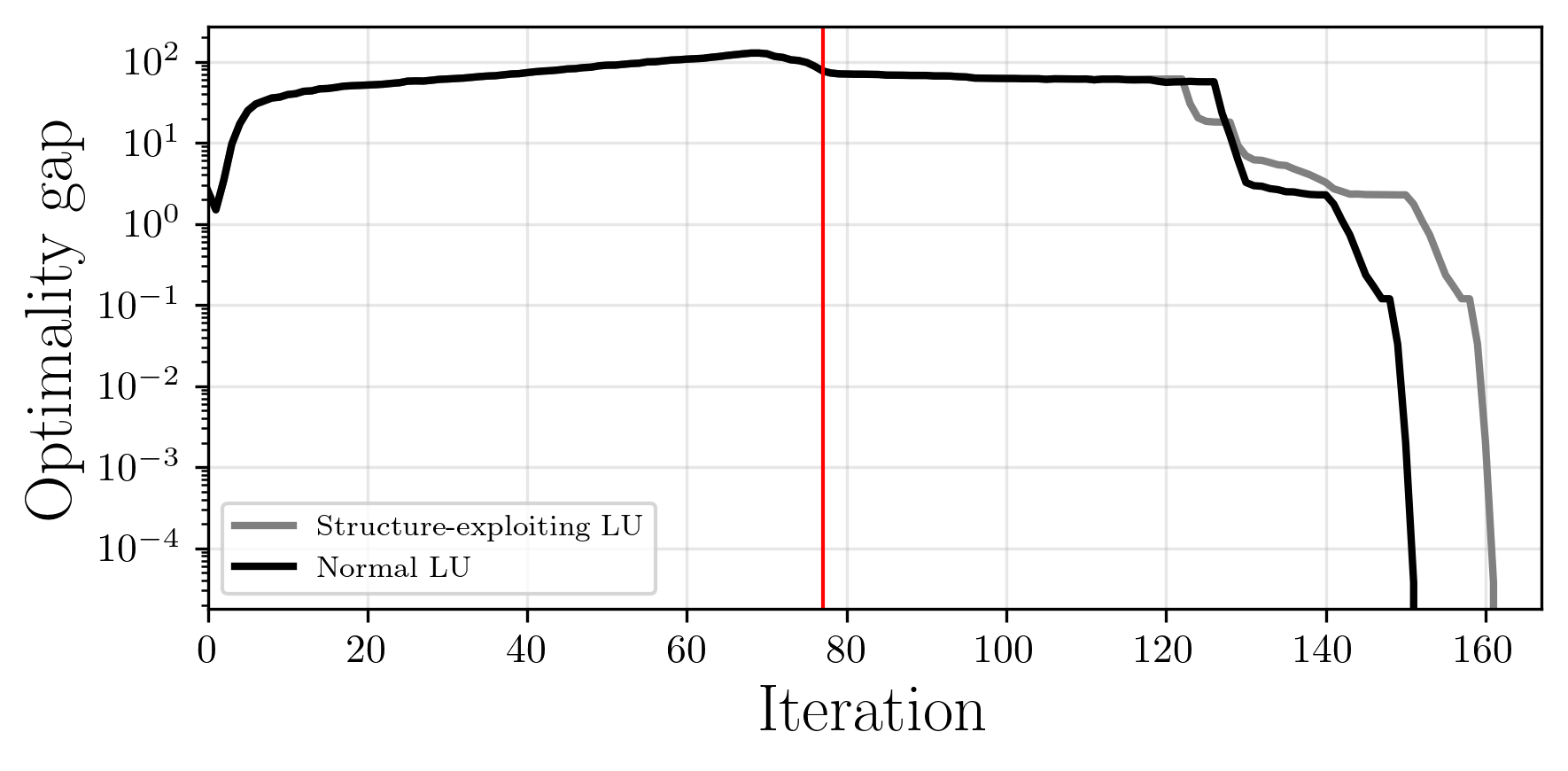}
    }

    \subfloat{
    \includegraphics[width=0.8\linewidth]{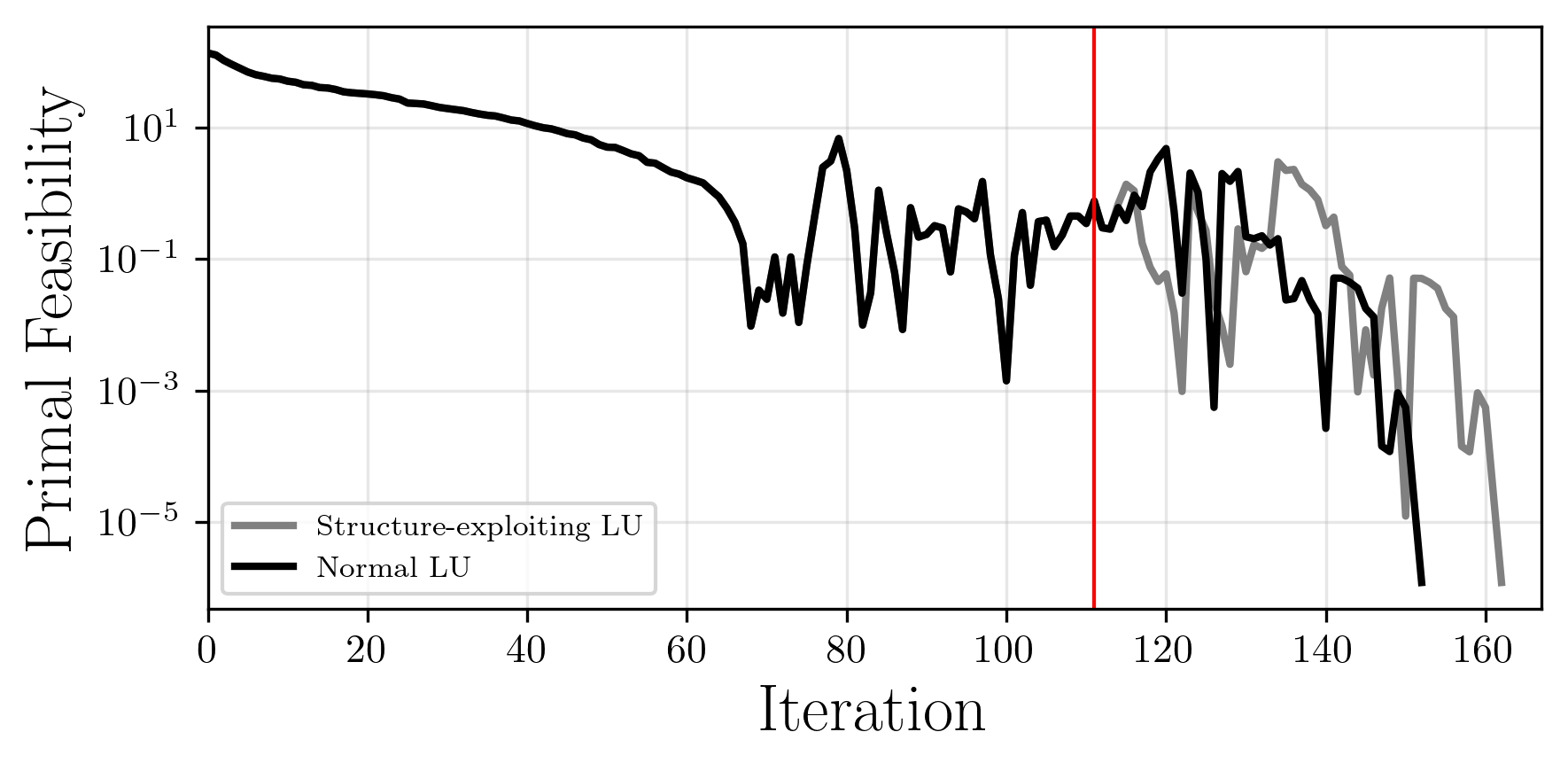}
    }

    \subfloat{
    \includegraphics[width=0.8\linewidth]{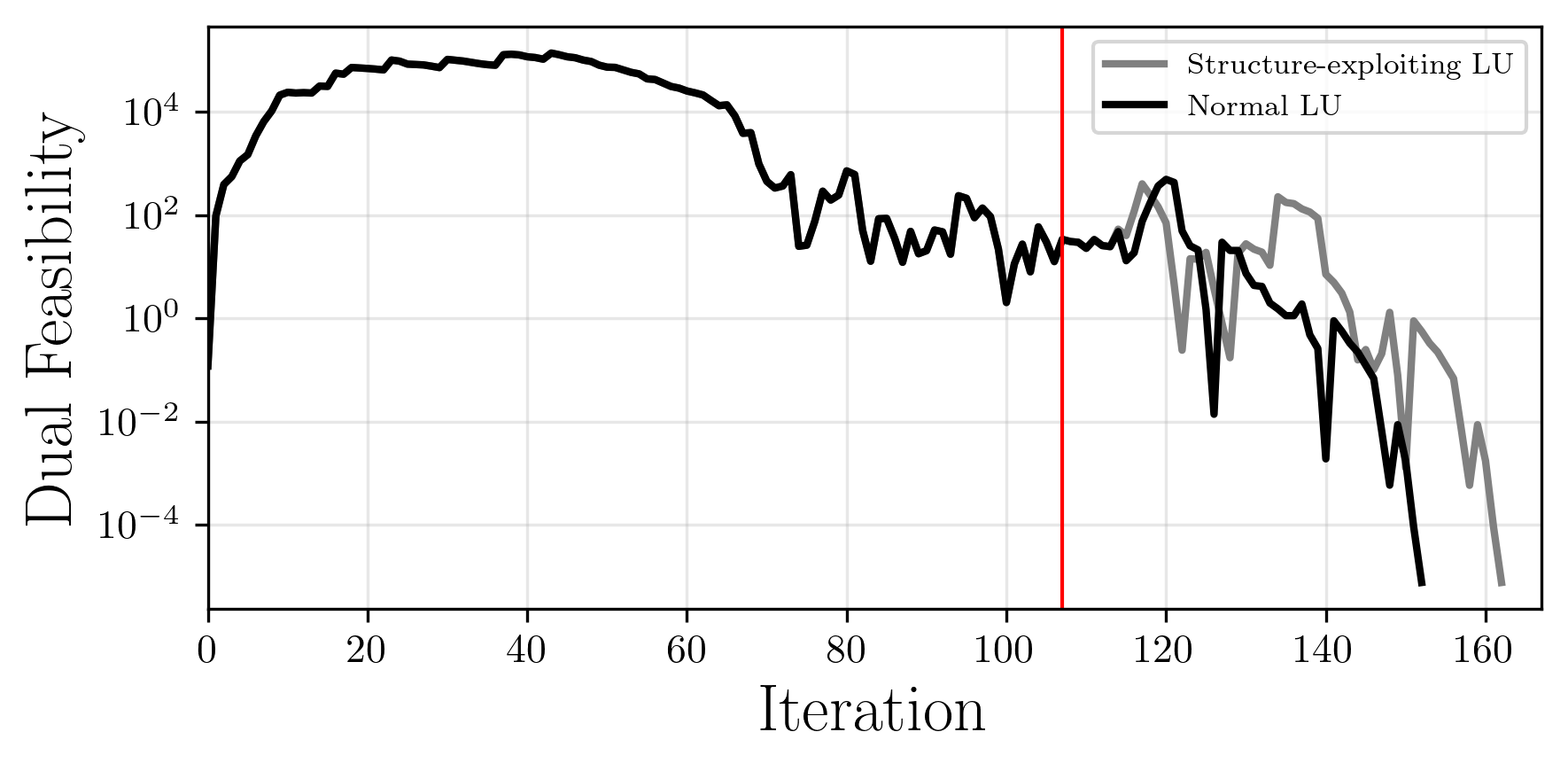}
    }

    \caption{Convergence plots of the walking gait OCP example. For each value plotted, a vertical red line indicates the first iteration at which a relative difference larger than $10^{-10}$ is observed. }
    \label{fig:convergence}
\end{figure}

For both OCP examples, Table \ref{tab:ocp_results} shows the number of solver iterations and computation times. For the 8-pendulum example, both the normal LU decomposition and the structure-exploiting LU decomposition lead to the same solver iterations, but the structure-exploiting LU decomposition leads to a 12\% decrease in the computation time of the search direction. For the walking gait example, the structure-exploiting LU decomposition leads to a reduction in computation time of the search direction per iteration of almost 20\%. However, both LU decompositions can select different pivots which can lead to numerical differences in iterations. Even though both approaches converge to the same solution, the structure-exploiting LU decomposition required more iterations for this example. Figure \ref{fig:convergence} shows the convergence behavior for the walking gait example. 



\begin{table}
    \centering
    \rowcolors{2}{gray!30}{white}
    \caption{Comparison between normal and structure-exploiting LU results for the OCP examples. The time required to find the search direction is denoted by $t_\mathrm{sd}$.}
    \label{tab:ocp_results}
    \subfloat{
    \begin{tabular}{l|ccc}
        \textbf{8-pendulum} & \begin{tabular}{@{}c@{}}Normal\\LU\end{tabular} & \begin{tabular}{@{}c@{}}Structure-exploiting\\LU\end{tabular} & \begin{tabular}{@{}c@{}}Relative \\ difference\end{tabular} \\
        \hline
        \# iterations & 119 & 119 & 0\% \\
        $t_\mathrm{total}$ (s) & 4.55 & 4.17 & -8.43\% \\
        $t_\mathrm{sd}$ (s) & 3.19 & 2.81 & -12.06\% \\
        $t_\mathrm{sd}$ per iter. (ms) & 27 & 24 & -12.06\%
    \end{tabular}
    }

    \rowcolors{2}{gray!30}{white}
    \subfloat{
    \begin{tabular}{l|ccc}
        \textbf{Walking gait} & \begin{tabular}{@{}c@{}}Normal\\LU\end{tabular} & \begin{tabular}{@{}c@{}}Structure-exploiting\\LU\end{tabular} & \begin{tabular}{@{}c@{}}Relative \\ difference\end{tabular} \\
        \hline
        \# iterations & 151 & 161 & 6.21\% \\
        $t_\mathrm{total}$ (s) & 20.49 & 19.72 & -3.9\% \\
        $t_\mathrm{sd}$ (s) & 11.34 & 10.09 & -12.33\% \\
        $t_\mathrm{sd}$ per iter. (ms) & 75.1 & 62.7 & -19.77\%
    \end{tabular}
    }
\end{table}

\section{Conclusions} \label{sec:conclusion}

This paper illustrates how structure-exploiting solvers such as Fatrop can benefit from exploiting not only the block-matrix structure of the KKT system but also the micro-structure within those block matrices. We show that reformulating general dynamics constraints to fit the problem structure of \eqref{eq:fatrop-ocp} introduces additional structure within the equality constraint Jacobian that can be exploited during the linear algebra computations.

The key contribution is the implementation of a structure-exploiting LU decomposition algorithm that exploits this structure and its integration within the Riccati recursion.

Numerical results demonstrate the effectiveness of the proposed method. On the 8-pendulum example, the structure-exploiting LU decomposition achieves an 8\% reduction in total computation time and a 12\% reduction in search direction computation time per iteration, while maintaining the same number of solver iterations. On the walking gait example, the speedup per iteration is more significant at nearly 20\%, though at the cost of additional solver iterations. The speedup is most pronounced for relative zero block sizes around 0.25.

The modified LU decomposition does not guarantee exactly the same solver iterations due to different pivot selections leading to numerical differences. The effect of the modified LU decomposition on the number of solver iterations should be further studied.

Another direction of interest is the automatic rearrangement of constraints and controls to optimize the relative size of the top-right zero block. This can be done if the sparsity structure of the equality constraint jacobian is known, which can easily be obtained when using a tool like CasADi\cite{casadi} to interface with the solver. Additionally, given the problem dimensions, automatic online selection of LU decomposition should be implemented such that the user can make abstraction of the underlying linear algebra computations while still benefiting from the speedup when possible.

\bibliography{bibliography.bib}

@misc{campd,
      title={Accelerated Multi-Modal Motion Planning Using Context-Conditioned Diffusion Models}, 
      author={Edward Sandra and Lander Vanroye and Dries Dirckx and Ruben Cartuyvels and Jan Swevers and Wilm Decré},
      year={2026},
      eprint={2510.14615},
      archivePrefix={arXiv},
      primaryClass={cs.RO},
      url={https://arxiv.org/abs/2510.14615}, 
	  howpublished = {arXiv preprint arXiv:2510.14615},
}

@article{ricatti_rao,
	title = {Application of {Interior}-{Point} {Methods} to {Model} {Predictive} {Control}},
	volume = {99},
	issn = {1573-2878},
	url = {https://doi.org/10.1023/A:1021711402723},
	doi = {10.1023/A:1021711402723},
	number = {3},
	journal = {Journal of Optimization Theory and Applications},
	author = {Rao, C. V. and Wright, S. J. and Rawlings, J. B.},
	month = dec,
	year = {1998},
	pages = {723--757},
}

@INPROCEEDINGS{ricatti-frison,
  author={Frison, Gianluca and Jørgensen, John Bagterp},
  booktitle={2013 IEEE International Conference on Control Applications (CCA)}, 
  title={Efficient implementation of the Riccati recursion for solving linear-quadratic control problems}, 
  year={2013},
  volume={},
  number={},
  pages={1117-1122},
  organization={},
  doi={10.1109/CCA.2013.6662901}
}

@article{hpipm,
title = {HPIPM: a high-performance quadratic programming framework for model predictive control⁎⁎This research was supported by the German Federal Ministry for Economic Affairs and Energy (BMWi) via eco4wind (0324125B) and DyConPV (0324166B), and by DFG via Research Unit FOR 2401.},
journal = {IFAC-PapersOnLine},
volume = {53},
number = {2},
pages = {6563-6569},
year = {2020},
note = {21st IFAC World Congress},
issn = {2405-8963},
doi = {https://doi.org/10.1016/j.ifacol.2020.12.073},
url = {https://www.sciencedirect.com/science/article/pii/S2405896320303293},
author = {Gianluca Frison and Moritz Diehl}
}

@article{blasfeo,
author = {Frison, Gianluca and Kouzoupis, Dimitris and Zanelli, Andrea and Diehl, Moritz},
year = {2017},
month = {04},
pages = {},
title = {BLASFEO: Basic Linear Algebra Subroutines for Embedded Optimization},
volume = {44},
journal = {ACM Transactions on Mathematical Software},
doi = {10.1145/3210754}
}

@article{ipopt,
	title = {On the implementation of an interior-point filter line-search algorithm for large-scale nonlinear programming},
	volume = {106},
	issn = {1436-4646},
	url = {https://doi.org/10.1007/s10107-004-0559-y},
	doi = {10.1007/s10107-004-0559-y},
	number = {1},
	journal = {Mathematical Programming},
	author = {Wächter, Andreas and Biegler, Lorenz T.},
	month = mar,
	year = {2006},
	pages = {25--57},
}

@article{qpoases,
	title = {{qpOASES}: a parametric active-set algorithm for quadratic programming},
	volume = {6},
	issn = {1867-2957},
	url = {https://doi.org/10.1007/s12532-014-0071-1},
	doi = {10.1007/s12532-014-0071-1},
	number = {4},
	journal = {Mathematical Programming Computation},
	author = {Ferreau, Hans Joachim and Kirches, Christian and Potschka, Andreas and Bock, Hans Georg and Diehl, Moritz},
	month = dec,
	year = {2014},
	pages = {327--363},
}

@article{fatrop-ricatti,
author = {Vanroye, Lander and De Schutter, Joris and Decré, Wilm},
title = {A generalization of the Riccati recursion for equality-constrained linear quadratic optimal control},
journal = {Optimal Control Applications and Methods},
volume = {45},
number = {1},
pages = {436-454},
doi = {https://doi.org/10.1002/oca.3064},
url = {https://onlinelibrary.wiley.com/doi/abs/10.1002/oca.3064},
eprint = {https://onlinelibrary.wiley.com/doi/pdf/10.1002/oca.3064},
year = {2024}
}

@INPROCEEDINGS{fatrop,
  author={Vanroye, Lander and Sathya, Ajay and De Schutter, Joris and Decré, Wilm},
  booktitle={2023 IEEE/RSJ International Conference on Intelligent Robots and Systems (IROS)}, 
  title={FATROP: A Fast Constrained Optimal Control Problem Solver for Robot Trajectory Optimization and Control}, 
  year={2023},
  volume={},
  number={},
  pages={10036-10043},
  organization={},
  doi={10.1109/IROS55552.2023.10342336}}

@Inbook{dea-multibody,
author="Arnold, Martin",
title="DAE Aspects of Multibody System Dynamics",
bookTitle="Surveys in Differential-Algebraic Equations IV",
year="2017",
publisher="Springer International Publishing",
address="Cham",
pages="41--106",
isbn="978-3-319-46618-7",
doi="10.1007/978-3-319-46618-7_2",
url="https://doi.org/10.1007/978-3-319-46618-7_2"
}

@article{baumgarte_stabilization_1972,
	title = {Stabilization of constraints and integrals of motion in dynamical systems},
	volume = {1},
	issn = {0045-7825},
	url = {https://www.sciencedirect.com/science/article/pii/0045782572900187},
	doi = {https://doi.org/10.1016/0045-7825(72)90018-7},
	number = {1},
	journal = {Computer Methods in Applied Mechanics and Engineering},
	author = {Baumgarte, J.},
	year = {1972},
	pages = {1--16},
}

@article {biomechanical,
	author = {D{\textquoteright}Hondt, Lars and Afschrift, Maarten and De Groote, Friedl},
	title = {Stochastic optimal control simulations of walking: potential and perspective},
	elocation-id = {2026.03.19.712839},
	year = {2026},
	doi = {10.64898/2026.03.19.712839},
	publisher = {Cold Spring Harbor Laboratory},
	URL = {https://www.biorxiv.org/content/early/2026/03/20/2026.03.19.712839},
	eprint = {https://www.biorxiv.org/content/early/2026/03/20/2026.03.19.712839.full.pdf},
	journal = {bioRxiv}
}

@Article{casadi,
  author = {Joel A E Andersson and Joris Gillis and Greg Horn
            and James B Rawlings and Moritz Diehl},
  title = {{CasADi} -- {A} software framework for nonlinear optimization
           and optimal control},
  journal = {Mathematical Programming Computation},
  volume = {11},
  number = {1},
  pages = {1--36},
  year = {2019},
  publisher = {Springer},
  doi = {10.1007/s12532-018-0139-4}
}
\nocite{*}

\appendix
\section{Verification of Structure-Exploiting LU decomposition} \label{app:structure-exploiting-lu}
To verify Algorithm \ref{alg:structure-exploiting-lu-carry}, compute both the left and right-hand side of the equation $P \overline{G}_u Q = L U$. The right-hand side is given by
\begin{align}
    \begin{split}
        \begin{bmatrix}
            K_1 & & \\
            K_7 & L_2 & \\ 
            K_2 &  & \bm{I}_{n_g-\rho_1}
        \end{bmatrix}
        \begin{bmatrix}
            V_1 & \tilde{V}_2 \\
            & U_2 \\
            & \bm{0}
        \end{bmatrix}
        = 
        \begin{bmatrix}
            & K_1 V_1 & K_1 \tilde{V}_2 \\
            & K_7 V_1 & K_7 \tilde{V}_2 + L_2 U_2 \\
            & K_2 V_1 & K_2 \tilde{V}_2
        \end{bmatrix}
    \end{split}
    \label{eq:right-hand-side-carry}
\end{align}
where
\begin{align}
    K_7 \tilde{V}_2 + L_2 U_2 
    &= P_2 \begin{bmatrix} K_3\\K_4\end{bmatrix} \begin{bmatrix} V_2 & \bm{0} \end{bmatrix} Q_2 + P_2 \begin{bmatrix} K_5 & M_3\\ K_6 & M_5 \end{bmatrix} Q_2\\
    &= P_2 \begin{bmatrix} K_3 V_2 + K_5 & M_3\\ K_4 V_2 + K_6 & M_5 \end{bmatrix} Q_2
    \label{eq:carry-verify}
\end{align}

The left-hand side can be computed using the definitions
\begin{align}
    X_1 &\coloneqq \begin{bmatrix} P_1 & \\ & \bm{I}_{n_f + n_c} \end{bmatrix}
    \overline{G}_u
    \begin{bmatrix} Q_1 & \\ & \bm{I}_{n_z}\end{bmatrix}
    = \begin{bmatrix} P_1 M_1 Q_1 & \\ M_2 Q_1 & M_3\\ M_4 Q_1 & M_5\end{bmatrix}
\end{align}
with $P_1 M_1 Q_1 = \begin{bmatrix} K_1 & \\ K_2 & \bm{I} \end{bmatrix} \begin{bmatrix}V_1 & V_2\\ & \bm{0} \end{bmatrix}$ and 
\begin{align}
    X_2 &\coloneqq \begin{bmatrix} \bm{I}_\rho & & \\ & & \bm{I}_{n_f + n_c} \\ & \bm{I}_{n_g-\rho} & \end{bmatrix} X_1
    = \begin{bmatrix} \begin{bmatrix} K_1 V_1 & K_1 V_2 \end{bmatrix} & \\ M_2 Q_1 & M_3\\ M_4 Q_1 & M_5\\ \begin{bmatrix} K_2 V_1 & K_2 V_2 \end{bmatrix} & \end{bmatrix}
\end{align}
where 
\begin{align}
    \begin{bmatrix} M_2 \\ M_4 \end{bmatrix} Q_1 &= \begin{bmatrix} K_3 & K_5\\ K_4 & K_6 \end{bmatrix} \begin{bmatrix} V_1^{-1} & -V_1^{-1} V_2\\ & \bm{I}_{n_u-\rho} \end{bmatrix}^{-1}
    = \begin{bmatrix} K_3 & K_5\\ K_4 & K_6 \end{bmatrix} \begin{bmatrix} V_1 & V_2\\ & \bm{I}_{n_u-\rho} \end{bmatrix}\\
    &= \begin{bmatrix} K_3 V_1 & K_3 V_2 + K_5\\ K_4 V_1 & K_4 V_2 + K_6 \end{bmatrix}
\end{align}
Such that we can compute
\begin{align}
    P \overline{G}_u Q &= \begin{bmatrix} \bm{I}_\rho & & \\ & P_2 & \\ & & \bm{I}_{n_g-\rho} \end{bmatrix}
    X_2
    \begin{bmatrix}\bm{I}_\rho & \\ & Q_2 \end{bmatrix}\\
    &= \begin{bmatrix} \begin{bmatrix}K_1 V_1 & K_1 V_2\end{bmatrix} & \\ P_2 \begin{bmatrix} K_3 V_1 & K_3 V_2 + K_5\\ K_4 V_1 & K_4 V_2 + K_6 \end{bmatrix} & P_2 \begin{bmatrix} M_3 \\ M_5 \end{bmatrix}\\ \begin{bmatrix}K_2 V_1 & K_2 V_2\end{bmatrix} & \end{bmatrix}
    \begin{bmatrix}\bm{I}_\rho & \\ & Q_2 \end{bmatrix}\\
    &= \begin{bmatrix} 
        K_1 V_1 & K_1 \tilde{V}_2\\ 
        K_7 V_1 & P_2 \begin{bmatrix} K_3 V_2 + K_5 & M_3\\ K_4 V_2 + K_6 & M_5\end{bmatrix} Q_2\\
        K_2 V_1 & K_2 \tilde{V}_2 \end{bmatrix}
    \label{eq:left-hand-side-carry}.
\end{align}
Substituting \eqref{eq:carry-verify} in \eqref{eq:right-hand-side-carry}, we conclude equality of \eqref{eq:right-hand-side-carry} and \eqref{eq:left-hand-side-carry} and conclude algorithm \ref{alg:structure-exploiting-lu-carry} holds.

\end{document}